\documentclass{amsart}
\usepackage{amssymb}

\newcommand{\beq}[1]{\begin{equation}\label{#1}}
\newcommand{\eeq}{\end{equation}}

\theoremstyle{definition}

\title{On recurrence in G-spaces}

\author{  Igor Protasov, Ksenia Protasova (Kyiv University)}

\begin{document}

UDC {2010 MSC}: 37A05, 22A15, 03E05

\begin{abstract} We introduce and analyze the following general concept of
recurrence. Let $G$ be a group and let $X$ be a G-space with the action $G\times X\longrightarrow X$, $(g,x)\longmapsto gx$. For a family $\mathfrak{F}$ of subset of $X$ and $A\in \mathfrak{F}$, we denote $\Delta_{\mathfrak{F}}(A)=\{g\in G: gB\subseteq A$ for some $B\in \mathfrak{F},  \  B\subseteq A\}$, and say that a subset $R$ of $G$ is $\mathfrak{F}$-recurrent if $R\bigcap \Delta_{\mathfrak{F}} (A)\neq\emptyset$ for each $A\in \mathfrak{F}$.
\end{abstract}
\maketitle

\large

Let $G$ be a group with the identity $e$  and let $X$ be a $G$-space, a set with the action
$G\times X\longrightarrow X$, $(g,x)\longmapsto gx$. If $X=G$  and $gx$ is the product of $g$ and $x$ then $X$ is called a left regular $G$-space.

Given a $G$-space $X$, a family $\mathfrak{F}$ of subset of $X$ and $A\in \mathfrak{F}$, we denote
\begin{eqnarray}
\nonumber \Delta_{\mathfrak{F}}(A)=\{g\in G: gB\subseteq A  \text{  for  \ some }  B\in\mathfrak{F}, B\subseteq A\}.
\end{eqnarray}

Clearly,  $e\in \Delta_{\mathfrak{F}}(A)$ and if $\mathfrak{F}$ is upward directed $(A\in \mathfrak{F}$, $A\subseteq C$ imply $C\in \mathfrak{F})$ and if $\mathfrak{F}$ is $G$-invariant $(A\in \mathfrak{F},$ $g\in G$  imply $gA\in \mathfrak{F}$)   then
$$\nonumber \Delta_{\mathfrak{F}}(A)=\{g\in G: gA\cap A \in \mathfrak{F}\}, \ \
\Delta_{\mathfrak{F}}(A)=(\Delta_{\mathfrak{F}}(A))^{-1}.$$

If $X$ is a left regular $G$-space and $\emptyset\notin  \mathfrak{F}$ then $\Delta_{\mathfrak{F}}(A)\subseteq A A^{-1}.$

For a $G$-space $X$ and a family $\mathfrak{F}$ of subsets of $X$, we say that a subset $R$ of $G$ is $\mathfrak{F}$-{\it recurrent} if $\Delta_{\mathfrak{F}}(A)\cap R \neq \emptyset$  for every $A\in \mathfrak{F}$. We denote by $\mathfrak{R}_{\mathfrak{F}}$  the filter on $G$ with the base $\cap\{\Delta_{\mathfrak{F}}(A):  A\in \mathfrak{F}^{\prime}\},$  where $\mathfrak{F}^{\prime}$ is a finite subfamily of $\mathfrak{F}$, and note that, for an ultrafilter $p$ on $G$,  $\mathfrak{R}_{\mathfrak{F}}\in p$ if and only if each member of $p$ is $\mathfrak{F}$-recurrent.

The notion of an $\mathfrak{F}$-recurrent subset is well-known in the case in which $G$ is an amenable group, $X$ is a left regular $G$-space  and $\mathfrak{F}=\{ A\subseteq X: \mu(A)>0$  for some left invariant Banach measure $\mu$ on $X\}$. See [1]  and [2] for historical background.

Now we endow $G$  with the discrete topology and identity the Stone-$\check{C}$ech compactification $\beta G$ of $G$ with the set of all ultrafilters on $G$. Then the family $\{\overline{A}: A\subseteq G\}$, where $\overline{A}=\{p\in \beta G: A\in p\}$,  forms a base for the topology of $\beta G$. Given a filter $\varphi$ on $G$, we denote $\overline{\varphi}=\cap\{\overline{A}: A\in\varphi\}$.

We use the standard extension [3] of the multiplication on $G$ to the semigroup multiplication on $\beta G$. We take two ultrafilters $p, q\in \beta G$, choose $P\in p$ and, for each $x\in P$, pick $Q_{x}\in q$. Then $\cup_{x\in P} x Q_{x}\in pq$ and the family of these subsets forms a base of the ultrafilter $pq$.

We recall [4] that a filter $\varphi$ on a group $G$  is {\it left topological}  if $\varphi$ is a base at the identity $e$ for some (uniquely at defined) left translation invariant (each left shift $x\longmapsto gx$ is continuous) topology on $G$. If $\varphi$ is left topological then $\overline{\varphi}$ is a subsemigroup of $\beta G$ [4]. If $G=X$ and a filter $\varphi$ is left topological then $\varphi=\mathfrak{R}_{\varphi}$.
\vskip 5pt

{\bf Proposition 1.} {\it For every $G$-space $X$ and any family $\mathfrak{F}$ of subsets of $X$, the filter $\mathfrak{R}_{\mathfrak{F}}$ is left topological. }
\vskip 5pt

{\it Proof}. By  [4], a filter $\varphi$ on a group $G$ is left topological if and only if, for every $\Phi\in\varphi$, there is $H\in\varphi$, $H\subseteq\Phi$  such that, for every $x\in H$, $xH_{x} \subseteq \Phi$ for some $H_{x}\in\varphi$.

We take an arbitrary $A\in \mathfrak{F}$, put $\Phi=\triangle_{\mathfrak{F}}(A)$ and, for each  $g\in\triangle_{\mathfrak{F}}(A)$, choose $B_{g}\in \mathfrak{F}$ such that $gB_{g}\in A$. Then
$g\triangle_{\mathfrak{F}}(B_{g})\subseteq \triangle_{\mathfrak{F}}(A)$ so put $H=\Phi$.

To conclude the proof, let $ \ \  A_{1}, \ldots,A_{n}\in \mathfrak{F} \ \  $.
We denote $ \  \  \Phi_{1} = \triangle_{\mathfrak{F}}(A_{1}),\ldots,
 \Phi_{n} = \triangle_{\mathfrak{F}}(A_{n}),$  $\Phi=\Phi_{1} \cap \ldots \cap \Phi _{n} .$
 We use the above paragraph, to choose $H_{1},\ldots ,H_{n}$ corresponding to $\Phi_{1}, \ldots,\Phi_{n}$ and put $H=H_{1}\cap\ldots\cap H_{n}$.    $ \ \  \ \ \Box$

 \vskip 5pt

Let $X$ be a $G$-space  and let $\mathfrak{F}$ be a family of subsets of $X$. We say that a family $\mathfrak{F}^{\prime}$ of subsets of $X$ is $\mathfrak{F}$-{\it disjoint} if $A\cap B\notin \mathfrak{F}$  for any distinct $A,B \in \mathfrak{F}^{\prime}$.

A family $\mathfrak{F}^{\prime}$ of subsets of $X$ is called $\mathfrak{F}$-{\it packing large} if,  for each   $A\in\mathfrak{F}^{\prime}$, any $\mathfrak{F}$-disjoint family of subsets of $X$  of the form $gA,$ $g\in G$ is finite.

We say that a subset $S$  of a group $G$  is a $ \ \triangle_{\omega}$-{\it set}  if  $ \ e\in A$ and every infinite subset $Y$ of $G$ contains two distinct elements $x,y$ such that $x^{-1} y\in S$ and $y^{-1} x\in S$.
\vskip 5pt

{\bf Proposition 2.} {\it Let $X$  be a $G$-space and let $\mathfrak{F}$ be a $G$-invariant upward directed family of subsets of $X$. Then $\mathfrak{F}$ is $\mathfrak{F}$-packing large if and only if, for each $A\in \mathfrak{F}$, the subset $ \ \triangle_{\mathfrak{F}}(A) \ $ of $ \ G$ is a $ \ \triangle_{\omega}$-set}.

\vskip 5pt

{\it Proof}. We assume that $\mathfrak{F}$   is   $\mathfrak{F}$-packing large  and take an arbitrary infinite subset $Y$ of $G$. Then we choose distinct $g,h\in Y$ such that $gA\cap hA\in \mathfrak{F}$, so $g^{-1} h\in \triangle_{\mathfrak{F}}(A)$, $ h g\in \triangle_{\mathfrak{F}}(A)$ and $\triangle_{\mathfrak{F}}(A)$ is a $\triangle_{\omega}$-set.

Now we suppose that $\triangle_{\mathfrak{F}}(A)$   is  a $\triangle_{\omega}$-set and take an arbitrary infinite subset $Y$ of $G$. Then there are distinct $g,h\in Y$ such that $g^{-1} h\in \triangle_{\mathfrak{F}}(A)$ so $g^{-1} hA \cap A \in \mathfrak{F}$ and $gA\cap hA \in \mathfrak{F}$. It follows that the family $\{gA: g\in Y\}$ is not $\mathfrak{F}$-disjoint. $ \ \  \ \Box$
\vskip 5pt

{\bf Proposition 3.} {\it For every infinite group $G$, the following statements hold
\vskip 5pt

$(i)$ a subset $A\subseteq G$ is a $\triangle_{\omega}$-set if and only  if $e\in A$ and every infinite subset $Y$ of $G$ contains an infinite subset $Z$ such that $x^{-1} y\in A$, $y^{-1} x\in A$ for any distinct $x,y\in Z$;

$(ii)$ the family $\varphi$ of all $\triangle_{\omega}$-sets of  $G$ is a filter;

$(iii)$ if $A\in\varphi$ then $G=FA$ for some finite subset $F$ of $G$.}

\vskip 5pt

{\it Proof}. $(i)$ We assume that $A$ is a $\triangle_{\omega}$-set and define a coloring $\chi$ of $[Y]^{2}$, $\chi: [Y]^{2}\longrightarrow\{0,1\}$ by the rule: $\chi(\{x,y\})=1$ if and only if $x^{-1} y\in A$, $y^{-1} x\in A$. By the Ramsey theorem, there is an infinite subset $Z$ of $Y$ such that $\chi$ is monochrome on $[Z]^{2}$. Since $A$ is  a $\triangle_{\omega}$-set $\chi(\{x,y\})=1$  for all $\{x,y\}\in [Z]^{2}$.
\vskip 5pt

$(ii)$ follows from  $(i)$.
\vskip 5pt

$(iii)$ We assume the contrary and choose an injective sequence $(x_{n})_{n\in\omega}$ in $G$  such that $ \ x_{n+1}\notin  x_{i} A \ $ for each $ \ \ i\in\{0,\ldots  , n\}$, and denote $Y=\{x_{n}: n\in\omega\}$. Then $x_{m} ^{-1} x_{n} \in A \ \ $ for every  $\  \  m, n $,  $ \ \ m<n$, so $A$ is not a $\triangle_{\omega}$-set.  $ \   \  \Box$
\vskip 5pt

\textbf{Proposition 4.} \textit{Let $G$ be a infinite group and let $\varphi$ denotes the filter of all $\triangle_{\omega}$-sets of $G$. Then $\overline{\varphi}$ is the smallest closed subset of $\beta G$ containing all ultrafilters on $G$  of the form $q^{-1} q$, $q\in \beta G$, $g^{-1} = \{A^{-1} : A\in q\}$. }
\vskip 5pt

{\it Proof}. We denote by $Q$ the smallest closed subset of $\beta G$ containing all $q^{-1} q $,   $ \ \ q\in \beta G$. It follows directly from the definition of the multiplication in $\beta G$ that $p\in Q$ if and only if either $p$ is principal and $p=e$ or, for each $P\in p$, there is an injective sequence $(x_{n})_{n\in\omega}$ in $G$ such that $x_{m} ^{-1}  x_{n} \in P$ for all $m<n$.

Applying Proposition $3(i)$, we conclude that $ \ q^{-1} q \in \overline{\varphi} \ $ for each $ \ q\in \beta G \ $ so $ \ Q\subseteq \overline{\varphi}$. On the other hand, if $p\notin \overline{\varphi}$  then there is $P\in p$ such that $G\setminus P$ is a $\triangle_{\omega}$-set.  By above paragraph, $p\notin Q$ so $\overline{\varphi}\subseteq Q$.     $\   \   \ \Box$
\vskip 5pt

Now let $G$ be an amenable group, $X$ be a left regular $G$-space and $\mathfrak{F}=\{A\in X: \mu(A)>0$  for some left invariant Banach measure $\mu$ on $G\}$. For combinatorial characterization  of $\mathfrak{F}$ see [6].
Clearly, $\mathfrak{F}$ is upward directed $G$-invariant and $\mathfrak{F}$-packing large. By  Proposition 2, $\overline{\varphi}\subseteq \overline{\mathfrak{R}_{\mathfrak{F}}}$. By  Proposition 4, $\overline{\mathfrak{R}_{\mathfrak{F}}}$ contains all ultrafilters of the form $q^{-1}q$, $q\in \beta G$, so we get Theorem 3.14 from [1].
\vskip 5pt

We suppose that a $G$-space $X$ is endowed with a $G$-invariant probability measure $\mu$ defined on some ring of subsets of $X$. Then the family $\mathfrak{F}\{ A\subseteq X: \mu(B)> 0$
for some $B\subseteq A\}$  is $\mathfrak{F}$-packing large.

In particular, we can take a compact group $X$, endow $X$ with the Haar measure, choose an arbitrary subgroup $G$ of $X$ and endow $G$ with the discrete topology.
\vskip 5pt

Another example: let a discrete group $G$ acts on a topological space $X$ so that, for each $g\in G$, the mapping $X\longrightarrow X$, $(g,x)\longmapsto gx$ is continuous. We take a point $x\in X$, denote by $\mathfrak{F}$ the filter of all neighborhoods of $x$, and recall that $x$ is {\it recurrent} if, for every $U\in \mathfrak{F}$, there exists $g\in G\backslash \{e\}$ such that $gx\in U$. Clearly, $x$ is a recurrent point if and only if $G\setminus\{e\}$ if a set of $\mathfrak{F}$-recurrence, so by Proposition 1, $x$ defines some non-discrete left translation invariant topology on $G$.
\vskip 5pt

{\bf Proposition 5.} {\it Let $G$ be a infinite group, $A$  be a $\triangle_{\omega}$-set of $G$ and let $\tau$ be a left translation invariant topology on $G$ with continuous inversion $x\longmapsto x^{-1}$ at the identity $e$. Then the closure $cl_{\tau}A$  is a neighborhood of $e$ in $\tau$.}

\vskip 5pt

\textit{Proof}. On the contrary, we suppose that $cl_{\tau}A$  is not a neighborhood of $e$, put $U=G\setminus cl_{\tau}A$. Then $U$ is open and  $e\in cl_{\tau}U$.

We take an arbitrary $x_{0}\in U$ and choose an open neighborhood $U_{0}$ of the identity such that $x_{0} U_{0}^{-1}\subseteq U$. Then we take $x_{1} \in U_{0}\cap U$ and choose an open neighborhood $U_{1}$ of $e$ such that $U_{1}\subseteq U_{0}$ and  $x_{1} U_{1}^{-1}\subseteq U$. We take $x_{2}\in U_{1}\cap U$ and choose an open neighborhood $U_{0}$ of $e$ such that $U_{2}\subseteq U_{1}$ and $x_{2} U_{2} ^{-1}\subseteq U$  and so on. After $\omega$ steps, we get a sequence $(x_{n})_{n\in \omega}$ in $G$ such that $x_{n} x_{m} ^{-1} \in U$ for all $n<m$. We denote $Y=\{x_{n} ^{-1}: n\in \omega\}$. Then $(x_{n} ^{-1})^{-1} x_{m} ^{-1}\in A$ for all $n<m$, so $A$ is not a $\triangle_{\omega}$-set. $\   \  \Box$

A subset $A$ of an infinite group $G$ is called a {\it $\triangle_{<\omega}$-set} if $e\in A$ and there exists a natural number $n$ such that every subset $Y$ of $G$, $\mid Y\mid=n$  contains two distinct $x, y\in Y$ such that $x^{-1} y \in A$, $y^{-1}x\in A$. These subsets were introduced in [5] under name thick subsets, but thick subsets are well-known in combinatorics with another meaning [3]: $A$ is thick if, for every finite subset $F$ of, there is $g\in A$ such that $Fg\subseteq A$.
The family $\psi$ of all $\triangle_{<\omega}$-sets of $G$ is a filter [5], clearly, $\psi\subseteq \varphi$. Every infinite group $G$ has a $\triangle_{\omega}$-set but not $\triangle_{<\omega}$-set $A$: it suffices to choose inductively a sequence $(X_{n})_{n\in\omega}$ of subsets of $G$, $\mid X_{n}\mid=n$ such that $\bigcup_{n\in\omega}X_{n}^{-1}X_{n}$  has no infinite subsets of the form $Y^{-1} Y$ and put $A=\{e\}\cup (G\backslash \bigcup_{n\in\omega}  X_{n}^{-1}X_{n})$, so $\psi\subset\varphi$.

By analogy with Propositions 3 and 4, we can prove
\vskip 5pt

{\bf Proposition 6.} {\it Let $G$ be an infinite group and let $\psi$ be the filter of all $\triangle_{<\omega}$-subsets of $G$. Then $p\in \overline{\psi}$ if and only if either $p$ is principal  and $p=e$ or, for every $A\in p$, there exists a sequence $(X_{n})_{n\in\omega}$ of subsets of $G$, $|X_{n}|= n+1$,  $X_{n}= \{x _{n0}, \ldots , x _{nn}\}$  such that  $x _{ni}^{-1} x _{n}j\in A$  for  all $i<j\leq n$.}\vskip 5pt

Let $A$ be a subset of a group $G$ such that $e\in A$,
$A=A^{-1}$. We consider the Cayley graph $\Gamma_{A}$ with the set of vertices $G$
and the set of edges $\{\{x,y\}: x^{-1}y\in A, \ \  x\neq y\}$.
We recall that a subset $S$ of vertices of a graph is {\it independent}
if any two distinct vertices from $S$ are not incident. Clearly, $A$ is a $\triangle_{\omega}$-set
if and only if any independent set in $\Gamma_{A}$ is finite, and $A$ is $\triangle_{\omega}$-set
if and only if  there exists a natural number $n$ such that any independent set $S$ is of size $|S|<n$.

{\bf Problem 1.} {\it  Characterize all infinite graphs with only finite independent set of vertices.}

{\bf Problem 2.} {\it  Given a natural number $n$, characterize all infinite graphs such that any independent set $S$ of vertices is of size $|S|<n$}.

In the context of this note, above problems are especially interesting in the case of Cayley graphs of groups.

\vspace{5 mm}
CONTACT INFORMATION

\vspace{5 mm}

I.~Protasov: \\
Department of Cybernetics  \\
         National Taras Shevchenko  University of Kiev \\
         Academic Glushkov St.4d  \\
         03680 Kiev, Ukraine \\ i.v.protasov@gmail.com

\medskip

K.~Protasova:\\
Department of Cybernetics  \\
         National Taras Shevchenko  University of Kiev \\
         Academic Glushkov St.4d  \\
         03680 Kiev, Ukraine \\ ksuha@freenet.com.ua

\end{document}